\documentstyle[11pt]{article}
\begin{document}
\title{{Bounds of Ideal Class Numbers of Real  Quadratic
  \\   Function Fields}
\thanks{Project Supported by the NNSFC (No.19771052)}
 \author{ \mbox{}\vspace{0.3 cm}{WANG Kunpeng
 \thanks{Current address: State Key Laboratory of Information Security,
Graduate School of Chinese Academy of sciences, Beijing, 100039 }
  $\;\;$ and $\; \;$ ZHANG Xianke }\\  \small \it{Tsinghua University,
  Department of Mathematical Sciences, Beijing 100084, P. R. China }
       \\
     \small \it{E-mail:\ kunpengwang@263.net, \ xianke@tsinghua.edu.cn}\\
     \small \it{Fax:\ 086-010-62781785}  } }
\date{}
\maketitle
\parindent 24pt
\baselineskip 20pt
\parskip 0pt

\noindent
 {\bf Abstract}.
 Theory of continued fractions of functions
  is used  to give  lower bound for class
 numbers $h(D)$ of general real quadratic function fields
 $K=k(\sqrt D)$ over $k={\bf F}_q(T)$.
  For five  series of real quadratic function fields $K$,
the bounds of $h(D)$ are given more explicitly, e.g.,
 if  $\ D=F^2+c,$ \mbox{}\hspace{0.1cm} then
$\ h(D)\geq \mbox{deg}F /\mbox{deg} P;$ \hspace{0.1cm}
if $D=(SG)^2+cS,\ $ then $\ h(D)\geq
\mbox{deg}S / \mbox{deg} P;\; $ if $D=(A^m+a)^2+A,\ $
 then  $\ h(D)\geq
\mbox{deg}A / \mbox{deg} P, \; $ where $P$ is irreducible
 polynomial splitting in $K,\; c\in {\bf F}_q$ .
 In addition, three types of quadratic
 function fields $K$ are found to have ideal class
 numbers bigger  than one.  \par \vskip 0.2cm

  \par \vskip 0.1cm

\noindent
 MR (2000) Subject Classification: 11R58; 11R29; 14H05.\\
 China Library Classification: O156.2
\par \vskip 0.6cm

\noindent
 {\Large I. Introduction and Main Results}
\vskip 0.4cm

Suppose that $k={\bf F}_q(T)$ is the rational function field in indeterminate
  (variable) $T$ over ${\bf F}_q$, the finite field with $q$
 elements ($q$ is a power of odd prime number).
  Let $R={\bf F}_q[T]$ be the polynomial ring of $T$ over
  ${\bf F}_q$. Any finite algebraic
  extension $K$ of $k$ is said to be an (algebraic)  function
  field.
  The integral closure of $R$ in
   $K$ is said to be the ring (domain) of
  integers of $K$, and is denoted by $ {{\cal O}}_K$,
 which
  is a Dedekind domain.
  The fractional ideals of $  {{\cal O}}_K$ form a multiplication group
  $  {{\cal I}}_K $. Let $ {{\cal P}}_K$ denotes
  the principal ideals in
  $  {{\cal I}}_K $. Then the quotient group
   $H(K)=  {{\cal I}}_K /  {{\cal P}}_K $
    is said to be the ideal
   class group of $K$.  And  $h(K)= \#H(K)$
   (the order of $H(K)$) is said to be the
   ideal class number of $K$. \vskip 2pt

 A quadratic extension of $k$ could be expressed as
  $K=k(\sqrt{D})$, where $D\in R$ is a polynomial
  which is not a square (we could also assume $D$ is
  square-free). If in addition $D$ is monic with even degree,
   then $K=k(\sqrt{D})$ is said to be  a real quadratic
    function field.
For real quadratic number fields $K$,  Mollin in 1987
 obtained a lower bound of ideal class numbers $h(D)=h(K)$
  by evaluating the fundamental unit (see [1]). Feng and Hu in [2]
   obtained a similar results for function fields;
   they also gave an explicit
   bound of $h(K)$
    for $K=k(\sqrt{F^2+c}) $, where $c\in {\bf F}_q^\times$.
There are also other works using continued fractions to study
quadratic number fields(see [3-6]).
 We here give  a  theorem on lower bound of
 $h(K)$ for general real quadratic function fields $K$,
 and  obtain explicit lower bound of $h(K)$ for six types of $K$
 including   the fields $K=k(\sqrt{F^2+c}) $.
 \vskip 0.2cm

 E. Artin in [7] began to use continued fractions to study
 quadratic function fields. In [8] we re-developed the
 theory of continued fractions of algebraic functions in
 an elementary and practicable way and studied some
properties of them, which will be  used here.   \vskip 0.2cm

 Suppose that $D$ is a monic square-free  polynomial with
degree
 $2d.$ By [8] we know $\sqrt{D}$ has an expansion of
 (simple) continued
fraction:  $\ \sqrt{D}=[a_0,a_1,\cdots].\ $
 Then $\alpha_i$ $=[a_i,a_{i+1},\cdots]$
 is said to be the   $i-$th complete quotient which
  could be expressed as
  $$\hskip 4cm \alpha_i=(\sqrt{D}+P_i)/Q_i\; \; \hskip 2cm
  (P_i,  Q_i \in R),\; \; $$
  and $Q_i$ is said to be the $i-$th complete denominator.
 The fraction $p_{i}/q_{i}=[a_0,a_1,\cdots,a_{i}]$
  is named the $i-$th
convergent. There is a positive integer $ \ell $
 such that
   $a_{n+\ell } = a_n  \; $ for any
 $1 \leq n\in {\bf Z}\ $ (The minimal $\ell $ having
 this property is called the $\underline{period}$
  of the continued
 fraction). So the continued fraction could be written as
 \par \vskip 0.1cm  \hskip 4.5cm
 $\ \sqrt{D}=[a_0, \overline{a_1, \cdots, a_{\ell}}], \ $
 \par \vskip 0.1cm
 \noindent
 where the underline part denotes a period, and we have
 also $ a_{\ell -i} = a_i $ for $ 0< i < \ell .\ $
Further more,
 there is a positive $v\in {\bf Z}$ such that
 $a_{n+v} = c a_n  \hskip 0.1cm  {\rm or}
 \hskip 0.1cm   c^{-1}  a_n\ $ for any
 $1 \leq n\in {\bf Z}$,
 where $c\in  {\bf F}_q^\times$.
 The minimal $v$ having
 this property is called the
 $\underline{{quasi-period}}$
  ($v$ is also the minimal integer$(>1)$ such that
  $Q_v \in {\bf F}_q^\times)$). We have
    $v= \ell/2\; {\rm or} \; \ell$
  (see [8]).  \vskip 0.2cm

\noindent
 {\bf Theorem 1.} Suppose that $K=k(\sqrt{D})$ is a real
quadratic function field, $\mbox{deg}D=2d$, and $P\in R$ is an
irreducible polynomial
 splitting in $K$. Then the ideal class number $h(D)$
of $K$ has a factor $h_1$ satisfying    \par \vskip 0.1cm
 \hskip 3cm $h_1= \mbox{deg} Q_i /\mbox{deg} P  \; \;  (1\leq i\leq v),
  \quad \mbox{or} \quad  h_1\geq d/\mbox{deg} P.$
   \par
\hskip -0.8cm In particular, we have
\par \hskip 2.5cm $
h(D)\geq {{min}}_{_{0< i <v}} \{\mbox{deg} Q_i /\mbox{deg}  P\ ,
\quad
  d/\mbox{deg}  P\}. $
     \par \vskip 0.5cm

\noindent
 {\bf Theorem 2.} Suppose that   $K=k(\sqrt{D})$ is a
real quadratic function field,
  $D\in R$ is square-free  with
$\mbox{deg}D=2d$, and $P\in R$ is an irreducible polynomial
 splitting in $K$. Then the ideal class number $h(D)$ of $K$ has the following
lower bound: \par  \vskip 0.1cm  \mbox{}\hspace{1cm} (1) If
$D=F^2+c,$ \mbox{}\hspace{1.06cm} then \quad
$h(D)\geq \mbox{deg}F
/\mbox{deg} P;$ \par  \vskip 0.1cm  \mbox{}\hspace{1cm}
 (2) If $D=(SG)^2+cS,$ \quad
then \quad $h(D)\geq \mbox{deg}S / \mbox{deg} P.$
\par \vskip 0.1cm \hskip -0.8cm where $ c \in {\bf F}_q^\times $,  $G\in R$
with $\mbox{deg}G\geq 1.$
\par \vskip 0.42cm

\noindent
 {\bf Theorem 3.}  Let  $ D\in R $ be square-free
polynomials
 as the followings, where $\ a\in R -{\bf F}_q,\ $
$  A=2a+1 $   is monic,   $ m $ is any positive integer.
Assume  $ P $ is an irreducible polynomial in $R$  splitting
 in $ K=k(\sqrt{D})$. Then the ideal class
number $h(D)$ of $K$ has  bound as the following: \par
   \vskip 0.1cm \mbox{}\hspace{1cm}(1) If $D=(A^m+a)^2+A,$
    \mbox{}\hspace{0.84cm} then \quad  $h(D)\geq
\mbox{deg}A / \mbox{deg} P;$\par
  \vskip 0.1cm \mbox{}\hspace{1cm}(2) If $D=(A^m-a)^2+A,$
\mbox{}\hspace{0.84cm} then \quad $h(D)\geq \mbox{deg}A /
 \mbox{deg} P;$\par
 \vskip 0.1cm\mbox{}\hspace{1cm}(3) If $D=(A^m+a+1)^2-A,$
 \quad then \quad $h(D)\geq \mbox{deg}A / \mbox{deg} P;$\par
  \vskip 0.1cm  \mbox{}\hspace{1cm}(4) If $D=(A^m-a-1)^2-A,$
 \quad then \quad $h(D)\geq \mbox{deg}A / \mbox{deg} P.$
 \vskip 0.3cm

 Now it is easy to find  fields $K=k(\sqrt{D})$ with
  class numbers $h(D)>1$. In the
  following Corollaries 1-3, we assume $P\in R$ is irreducible,
   $\  c\in {\bf F}_q^\times,\ $  and $\ (F/P)\ $ is the
   quadratic-residue symbol
 (i.e., $(F/P)=1$ or $-1$ according to $F$ is
  a quadratic residue modulo $P$ or not).  \vskip 0.3cm

\noindent
 {\bf Corollary 1.} Let $\ D=(PG)^2+c,\ $ $\
\mbox{deg}G\geq 2,\ $ $\ (c/P)=1, \quad $  then \par \vskip 0.1cm
\hskip 4cm $h(D)\geq \mbox{deg}(GP)/\mbox{deg}P>1.$
\par \vskip 0.26cm

\noindent
 {\bf Corollary 2.} Let
 $\ D=(SHP)^2+cS,\ $  $\ \mbox{deg}(S) > \mbox{deg}(P),\ $ $(cS/P)=1,\ $
  then \par \vskip 0.1cm   \hskip 4cm
   $ h(D)\geq \mbox{deg}(S)/\mbox{deg}P>1$.
\par \vskip 0.26cm

\noindent
 {\bf Corollary 3.} Let $D$ be as in Theorem 3 with $\
A=SP,\ $
 $\ \mbox{deg}S\geq 2, \quad$
 then \par \vskip 0.1cm   \hskip 4cm
  $h(D)\geq \mbox{deg}(SP)/\mbox{deg}P>1.$
   \vskip 0.26cm

 As an example of Corollary 2, we have
 \par \vskip 0.1cm   \hskip 4cm
  $h(\ (T(T^m+1))^2+(T^m+1) \ )\ \geq \ m .$
  \par \vskip 0.1cm
 \hskip -0.8cm (i.e., we take $ P=T, \;$  $S=(T^m+1), \;$  $H=1,\;$
 $c=1$)
  \par \vskip 0.56cm

\noindent
 {\Large  II. Lemmas and Proofs of Theorems}
\par \vskip 0.4cm

 First, consider the expansion and property of continued
 fractions of  $\sqrt D,\  $  where  $ D\in R $ is a
 square-free  monic polynomial with even degree
 $\ \mbox{deg} (D)=2d. \ $  By [8] we know  there exist
uniquely determined  $\; f, r \in R $
 such that $\; D=f^2+r, \; $ and $f$ is monic,
  $\ \mbox{deg}f=d,\; $ $\mbox{deg}r<d. \;  $
\par \vskip 0.36cm

 The following process
 produces the expansion of simple continued fraction
 $\ \sqrt{D}=[a_0,a_1,\cdots] :$  \par  \vskip 2pt

1. Denote $D=f^2+r$ as above.  Put $a_0=f$, then
 $\sqrt{D}=a_0+\sqrt{D}-a_0,$
   thus $\alpha_1=1/(\sqrt{D}-a_0)=
   (\sqrt{D}+a_0)/({D}-a_0^2)=(\sqrt{D}+P_1)/Q_1$,
    where $P_1=a_0$, $Q_1=D-a_0^2$.  \vskip 0.2cm

2. Now  $\alpha_1=(f+P_1+\sqrt{D}-f)/Q_1.\ $ Assume
 $\ f+P_1=a_1 Q_1+r_1,\; \; \mbox{deg}\; r_1< \mbox{deg} Q_1.\ $
  Then $\ \alpha_1 =a_1+(\sqrt{D}-(f-r_1))/Q_1.\ $
 Thus $\ \alpha_2$ $=Q_1/(\sqrt{D}-(f-r_1))$
$=(\sqrt{D}+P_2)/Q_2,\ $
 where $\ P_2=(f-r_1),\ $
 $\ Q_2=({D}-(f-r_1)^2)/Q_1=(D-P_2^2)/Q_1.\ $
   We see  $\ P_2\in R.\ $ Since $P_2=f-r_1=a_1 Q_1-P_1,\ $ so
   $\ D-P_2^2\equiv D-P_1^2\equiv 0\;  (\mbox{mod}Q_1),\; $
   thus $\; Q_2\in R.$
 \vskip 0.1cm

 Proceed continually, we could obtain
 the simple continued fraction of
  $\sqrt{D}$
 (see [8]).  \vskip 0.3cm

\noindent
 {\bf  Lemma 1}$^{[8]}$. The Diophantine equation $\
X^2-DY^2=G\ $  has a primary solution if and only if
  $ \ G=(-1)^i Q_i \ $   for some
   $\ 0\leq i \leq \ell ,\ $ where  $ Q_i $   is the
  $i-$th complete denominator of
  the continued fraction of   $ \sqrt{D}, $
   $D\in R $ is a monic
 square-free polynomial with even degree,
 $\ G\in R \ $ and
  $ \ \mbox{deg}  G < \frac{1}{2}\mbox{deg} D\ $.
(A solution $(X, Y)$  is primary
if $(X, Y)=1,\ X,Y\in R $ ). \par \vskip 0.2cm

\noindent
 {\bf Proof of Theorem 1 }.   Assume
  $\ (P)=\wp \overline\wp,\ $   where
 $\ \wp\not=\overline\wp\ $  are prime ideals of  $K$.  Let
 $ h=h(D) $   be the ideal class number of   $ K $,
 then   $\ \wp^h\ $ is a principal ideal.  Suppose that
   $\ m \leq h\ $  is the minimal positive integer such that
    $\ \wp^m \ $   is principal,  then  $\ m\ $  is a factor
    of $h$.
 Since $\ \{1,\ \sqrt D \}\ $ is an  integral basis for $K$,
so we may assume
 $\ \wp^m=(U+V\sqrt{D})\ $  with  $\  U,V \in R .\ $
Taking norm on both sides, we obtain an equation of ideals of
  $ k $ :  $\  ( U^2-DV^2 )=( P^m ).\ $
 So  $\ \ U^2-DV^2=cP^m \ $
 $\ (c\in {\bf F}_q^\times)\ $ since the unit group of $k$
 (or $R$) is just ${\bf F}_q^\times .\ $\par
 We assert that $U$ and $V$ must be relatively prime;
 otherwise, if $\ (U,\ V)=C\in R\ $
  is not a constant, put $\ U_1=U/C,\ V_1=V/C,\ $  then
  $\ \wp^m=(U+V\sqrt{D})=(C)(U_1+V_1\sqrt{D})\ $,
 by the uniqueness of factorization of ideals,  we must have
    $\ (U_1+V_1\sqrt{D})=\wp^n \ $
  for some $ n<m ,\ $
 which  contradicts  to the minimal assumption of  $ m $.\par
  Thus we know
$\ (U, V)\ $  is a primary solution of $\ X^2-DY^2=cP^m \ $.
First assume $ \mbox{deg} P^m <d ,\ $
 then by Lemma 1 we know that
 $ \; cP^m=(-1)^i Q_i ,\; $
   $ \ \mbox{deg} P^m=\mbox{deg} Q_i \ $
for some $i$ with $\; 0\leq i \leq \ell. \; $
   Thus  we have $ \ m =\mbox{deg} Q_i /\mbox{deg}  P\ $
   for some  $\; 0\leq i \leq v \; $ by the definition of
   quasi-period $\ v.\ $
Secondly, assume  $ \mbox{deg} P^m \geq d $, then we have
directly   $ \ m \geq    d /\mbox{deg}  P.\ $
\par \vskip 0.4cm

\noindent
 {\bf Proof of Theorem 2}. (1) It is easy to get the
expansion of simple continued fraction
 $ \sqrt{D}=\sqrt{F^2+c}=[F,\ \overline{2F/c,\ 2F}],\ $
and obtain the set of complete denominators:
   $\ (Q_0,\ Q_1,\ Q_2) = (1,\ c,\ 1).\ $
 Thus by Theorem 1 we know that  $h(D)$ has a positive
 factor $ h_1 $ satisfy
    $\ h_1\geq d / \mbox{deg}  P $, so
   $\ h(D)\geq \mbox{deg} F / \mbox{deg}  P $.\par
 \vskip 0.2cm

(2) Expand $ \sqrt{D}=\sqrt{ (SG)^2+cS}\ $ as simple continued
fraction: $\ \sqrt{D}=[SG,\ \overline{2G/c,\ 2SG}]. \  $  Its
period is
 $ \ell=2 $,  complete denominators are just
 $ (Q_0,\ Q_1, Q_2) = (1, \ cS,\ 1). $
 Thus by Theorem 1 we know   $\  h(D)\geq \mbox{deg} S /
\mbox{deg}  P $ (Note that $d={\mbox{deg}} D\geq \mbox{deg} S$
now).
\par \vskip 0.4cm

{\bf Proof of  Theorem 3}.
(1) The polynomial $D=(A^m+a)^2+A$ in the theorem
has good property which enables us to  expand $\ \sqrt{D} $
as a simple continued fraction
 $ \sqrt{D}=[a_0, \overline{a_1, \cdots,  a_\ell}] $.
 By Theorem 1, we need only to know a quasi-period of
  the expansion, i.e., $\ [a_0,\ \cdots,\ a_v].\ $
It turns out that this quasi-period is quite long and
demonstrates rules in three sections, so we will write down
 it in three sections and list    $\ a_n $,
  $P_n,\ Q_n\ (0\leq n \leq v)$. We need to distinguish
   four cases  \ $m=4t-2,\ 4t-1,\ 4t,\ 4t+1$ .
     \vskip 0.2cm

     The first section ($n=0, 1$) and the second section
     ($2\leq n \leq  4t+1$):
     \hskip 0.2cm \quad ($ 1\leq j\leq t $ )\\  \smallskip
\begin{tabular}{c|cccccccc}\hline
 $ n     $  &  $   0       $  &  $   1          $  &
 $  \cdots  $  &  $  4j-2          $  &
 $  4j-1        $  &  $  4j           $  &  $  4j+1
  $  &  $  \cdots $ \\ \hline
 $ P_n   $  &  $   0       $  &  $  A^m+a       $  &
  $  \cdots  $  &  $  A^m+a+1       $  &
 $  A^m-a-1     $  &  $  A^m+a+1      $  &  $  A^m-a-1
 $  &  $  \cdots $ \\ \hline
 $ Q_n   $  &  $   1       $  &  $   A          $  &  $
 \cdots  $  &  $  -2A^{m-2j+1}  $  &
 $  -A^{2j}     $  &  $  2A^{m-{2j}}  $  &  $  A^{2j+1}
 $  &  $  \cdots $ \\ \hline
 $ a_n   $  &  $   A^m+a   $  &  $  2A^{m-1}+1  $  &
  $  \cdots  $  &  $  -A^{2j-1}     $  &
 $  -2A^{m-2j}  $  &  $  A^{2j}       $  &  $  2A^{m-2j-1}
  $  &  $ \cdots $ \\ \hline
\end{tabular}

\par \vskip 0.26cm
 The third  section  is given in four cases
according to $m$: \\ \smallskip  % ~{\CC {\char35}:} \CCA
 \qquad  (i) For  $ m=4t-2:  $  \\ \begin{tabular}{c|cc}\hline
 $ n     $  &  $    4t+2           $  &  $      4t+3       $ \\ \hline
 $ P_n   $  &  $  A^m+a+1          $  &  $      A^m+a      $ \\ \hline
 $ Q_n   $  &  $     -2A           $  &  $  - 1/2    $ \\ \hline
 $ a_n   $  &  $  A^m- 1/2   $  &  $  -4A^m-4a       $ \\ \hline
\end{tabular}\   \par \vskip 0.2cm

(ii) For  $ m=4t-1 $ : \\ \begin{tabular}{c|ccccc}\hline
 $ n     $  &  $    4t+2      $  &  $      4t+3     $  &  $      4t+4    $  &  $    4t+5
    $  \\ \hline
 $ P_n   $  &  $   A^m+a+1    $  &  $  A^m-a-1      $  &  $    A^m+a+1   $  &  $    A^m+a
   $   \\ \hline
 $ Q_n   $  &  $ -2A^{m-2t-1} $  &  $  -A^{2t+2}    $  &  $    -A        $  &  $    -1
 $    \\ \hline
 $ a_n   $  &  $  -A^{m-2}    $  &  $  -2A          $  &  $  -2A^{m-1}+1 $  &  $
-2A^m-2a $     \\ \hline
\end{tabular}\  \par \vskip 0.2cm

 (iii) For  $ m=4t $:  \\ \begin{tabular}{c|cc}\hline
 $ n     $  &  $    4t+2                $  &  $      4t+3       $ \\ \hline
 $ P_n   $  &  $  A^m+a+1               $  &  $      A^m+a      $ \\ \hline
 $ Q_n   $  &  $  2A^{2-1}              $  &  $  - 1/2    $ \\ \hline
 $ a_n   $  &  $  -A^{m-1}- 1/2   $  &  $  -4A^m-4a       $ \\ \hline
 \end{tabular}  \   \par \vskip 0.2cm

(iv) For  $ m=4t+1 $:  \\ \begin{tabular}{c|ccccc}\hline
 $ n     $  &  $    4t+2      $  &  $      4t+3     $  &  $    4t+4       $  &  $    4t+5
        $ \\ \hline
 $ P_n   $  &  $   A^m+a+1    $  &  $  A^m-a-1      $  &  $  A^m+a+1      $  &  $    A^m+a
        $ \\ \hline
 $ Q_n   $  &  $ -2A^{m-2t-1} $  &  $   -A^{2t+2}   $  &  $  2A           $  &  $
 1/2   $ \\ \hline
 $ a_n   $  &  $  -A^{m-2}    $  &  $  -2A          $  &  $  -2A^{m-1}+1  $  &  $
-2A^m-2a       $ \\ \hline
\end{tabular}\  \par \vskip 0.3cm

(2) Similarly to (1), a quasi-period of the simple continued
fraction of $ \sqrt{D} $ is given here. \par
The first and second sections are combined given as
 the following: %~{\CC {\char35}:\CCA }
 \hskip 1cm ($ 1\leq j\leq t $   ) \\
\begin{tabular}{c|cccccccc}\hline
 $ n     $  &  $   0     $  &  $   1          $  &  $  \cdots  $  &  $  4j-2          $  &  $
4j-1        $  &  $  4j           $  &  $  4j+1         $  &  $  \cdots  $ \\ \hline
 $ P_n   $  &  $   0     $  &  $  A^m-a       $  &  $  \cdots  $  &  $  A^m-a-1       $  &  $
A^m+a+1     $  &  $  A^m-a-1      $  &  $  A^m+a+1      $  &  $  \cdots  $ \\ \hline
 $ Q_n   $  &  $   1     $  &  $   A          $  &  $  \cdots  $  &  $  2A^{m-2j+1}   $  &  $
-A^{2j}     $  &  $  -2A^{m-{2j}} $  &  $  A^{2j+1}     $  &  $  \cdots  $ \\ \hline
 $ a_n   $  &  $  A^m-a  $  &  $  2A^{m-1}-1  $  &  $  \cdots  $  &  $  A^{2j-1}      $  &  $
-2A^{m-2j}  $  &  $  -A^{2j}      $  &  $  2A^{m-2j-1}  $  &  $  \cdots $ \\ \hline
\end{tabular}\   \\

\par \vskip 0.16cm
The third section is given in four cases:

 (i) For  $ m=4t-2 $ : \\ \begin{tabular}{c|cc}\hline
 $ n     $  &  $    4t+2      $  &  $      4t+3       $ \\ \hline
 $ P_n   $  &  $  A^m-a-1     $  &  $      A^m-a      $ \\ \hline
 $ Q_n   $  &  $    2A        $  &  $     1/2   $ \\ \hline
 $ a_n   $  &  $  A^{m-1}-1   $  &  $     4A^m-4a     $ \\ \hline
\end{tabular}\\
\par \vskip 0.1cm

(ii) For   $ m=4t-1 $: \\ \begin{tabular}{c|ccccc}\hline
 $ n     $  &  $    4t+2      $  &  $      4t+3     $  &  $    4t+4      $  &  $    4t+5
       $ \\ \hline
 $ P_n   $  &  $   A^m-a-1    $  &  $     A^m+a+1   $  &  $  A^m-a-1     $  &  $    A^m-a
       $ \\ \hline
 $ Q_n   $  &  $  2A^{m-2t-1} $  &  $     -A^{2t}   $  &  $  -2A^{m-2t}  $  &  $
- 1/2    $ \\ \hline
 $ a_n   $  &  $  A^{m-2}     $  &  $  -2A          $  &  $  -A^{m-1}+1  $  &  $  -4A^m+4a
       $ \\ \hline
\end{tabular}\   \\
\par \vskip 0.16cm

(iii) For  $ m=4t $: \\ \begin{tabular}{c|cc}\hline
 $ n     $  &  $    4t+2      $  &  $      4t+3       $ \\ \hline
 $ P_n   $  &  $  A^m-a-1     $  &  $      A^m-a      $ \\ \hline
 $ Q_n   $  &  $    2A        $  &  $      1/2  $ \\ \hline
 $ a_n   $  &  $  2A^{m-1}-1  $  &  $     4A^m-4a     $ \\ \hline
\end{tabular}  \\
\par \vskip 0.16cm

(iv) For  $ m=4t+1 $: \\ \begin{tabular}{c|ccccc}\hline
 $ n     $  &  $    4t+2      $  &  $      4t+3     $  &  $    4t+4                $  &
 $    4t+5          $ \\ \hline
 $ P_n   $  &  $   A^m-a-1    $  &  $  A^m+a+1      $  &  $  A^m-a-1               $  &
 $    A^m-a         $ \\ \hline
 $ Q_n   $  &  $ 2A^{m-2t-1}  $  &  $   -A^{2t+2}   $  &  $  -2A                   $  &
 $    - 1/2   $ \\ \hline
 $ a_n   $  &  $  A^{m-2}     $  &  $  -2A          $  &  $  -A^{m-1}+ 1/2   $  &
 $  -4A^m+4a $ \\ \hline
\end{tabular}. \\  \vskip 0.3cm

(3) A quasi-period of the simple continued fraction of
 $ \sqrt{D} $ is given here in three sections similarly as in (1).
 The first and second
 sections of it are combined given as the following
  \hskip 1cm ($ 1\leq j\leq t $ ):\\
\begin{tabular}{c|cccccccc}\hline
 $ n     $  &  $   0        $  &  $   1           $  &  $  \cdots  $  &  $     4j-2
 $  &  $  4j-1       $  &  $   4j          $  &  $  4j+1        $  &  $ \cdots   $ \\ \hline
 $ P_n   $  &  $   0        $  &  $  A^m+a+1      $  &  $  \cdots  $  &  $   A^m+a
 $  &  $  A^m-a      $  &  $    A^m+a      $  &  $  A^m-a       $  &  $ \cdots   $ \\ \hline
 $ Q_n   $  &  $   1        $  &  $   -A          $  &  $  \cdots  $  &  $  -2A^{m-2j+1}
 $  &  $  -A^{2j}    $  &  $  -2A^{m-2j}   $  &  $  -A^{2j+1}   $  &  $ \cdots   $ \\ \hline
 $ a_n   $  &  $  A^m+a+1   $  &  $  -2A^{m-1}-1  $  &  $  \cdots  $  &  $  -A^{2j-1}
 $  &  $  -2A^{m-2j}  $  &  $  -A^{2j}      $  &  $  -2A^{m-2j-1} $ & $ \cdots $ \\ \hline
\end{tabular} \\

The third section is given in two cases.

(i) For  $ m=4t+1 $   or  $ m=4t+3 $: \\
\begin{tabular}{c|cc}\hline
 $ n     $  &  $    4t+2                $  &  $      4t+3        $ \\ \hline
 $ P_n   $  &  $    A^m+a               $  &  $      A^m+a+1     $ \\ \hline
 $ Q_n   $  &  $    -2A                 $  &  $     1/2    $ \\ \hline
 $ a_n   $  &  $  -A^{m-1}- 1/2   $  &  $    4A^m+4a+1     $ \\ \hline
\end{tabular}  \\

(ii) For  $ m=4t $   or  $ m=4t+2 $: \\ \begin{tabular}{c|ccccc}\hline
 $ n     $  &  $    4t+1     $  &  $     4t+2               $  &  $    4t+3         $ \\
\hline  $ P_n   $  &
  $   A^m-a     $  &  $    A^m+a               $  &  $  A^m+a+1        $ \\  \hline
 $ Q_n   $  &  $
-A^{2t+1}  $  &  $    -2A                 $  &  $  - 1/2    $ \\  \hline  $ a_n
   $  &  $  -2A  $  &  $
-A^{m-1}- 1/2   $  &  $  -4A^m-4a-4     $ \\  \hline
\end{tabular}  \\ \vskip 0.3cm

(4) We give a quasi-period of the simple continued fraction of
$ \sqrt{D} $ similarly as in (1). The first and second section are combined
given as the following:
\hskip 1cm ( $ 1\leq j\leq t $ )\\
\begin{tabular}{c|cccccccc}\hline
 $ n     $  &  $     0     $  &  $   1       $  &  $  \cdots  $  &  $     4j-2      $  &  $
4j-1        $  &  $    4j        $  &  $  4j+1         $  &  $ \cdots  $ \\ \hline
 $ P_n   $  &  $     0     $  &  $  A^m-a-1  $  &  $  \cdots   $  &  $   A^m-a       $  &  $
A^m+a       $  &  $    A^m-a     $  &  $  A^m+a        $  &  $ \cdots  $ \\ \hline
 $ Q_n   $  &  $     1     $  &  $   -A      $  &  $  \cdots  $  &  $  2A^{m-2j+1}  $  &  $
-A^{2j}     $  &  $   2A^{m-2j}  $  &  $  -A^{m-2j-1}  $  &  $ \cdots  $ \\ \hline
 $ a_n   $  &  $  A^m-a-1  $  &  $  -2A^m+1  $  &  $  \cdots  $  &  $   A^{2j-1}    $  &  $
-2A^{m-2j}  $  &  $   A^{2j}     $  &  $  -2A^{m-2j-1} $ &  $ \cdots $ \\ \hline
\end{tabular} \\

The third section is given in two cases.

 (i) For  $ m=4t+1 $ or
$ m=4t+3 $: \\ \begin{tabular}{c|cc}\hline
 $ n     $  &  $    4t+2               $  &  $    4t+3         $ \\  \hline  $ P_n   $  &
 $    A^m-a  $  &  $
A^m-a-1      $ \\  \hline  $ Q_n   $  &  $    2A                 $  &  $
- 1/2    $ \\  \hline
 $ a_n   $  &  $  A^{m-1}- 1/2   $  &  $  -4A^m+4a+4  $ \\  \hline
\end{tabular}  \\

(ii) For  $ m=4t $    or   $ m=4t+2 $: \\ \begin{tabular}{c|cc}\hline
 $ n     $  &  $    4t+2                $  &  $      4t+3        $ \\  \hline  $ P_n   $
&  $    A^m-a  $  &
 $      A^m-a-1     $ \\  \hline  $ Q_n   $  &  $    2A                  $  &  $
- 1/2   $ \\ \hline  $ a_n   $  &  $  A^{m-1}- 1/2    $  &  $  -4A^m+4a+4      $ \\  \hline
\end{tabular}. \\ \vskip 0.1cm

Consider the above simple continued fractions of $\sqrt D$ for the
four types of $D$, and check  the complete denominators $\{ {Q_n\}
(0 < n < v)} \ $  in a quasi-period, we find that the complete
denominator having the minimal degree is
  $ \pm 2A $ in all the cases.
Since deg$A<\frac{1}{2}$deg $D$, by Theorem 1 we know
 $ h(D)\geq \mbox{deg} A / \mbox{deg}  P. $
 This  proves Theorem 3.
 \vskip 24pt

\noindent
 {\large \bf References}

 \baselineskip 0pt
\parskip 0pt
\begin{description}
\item[[1]] R.A.Mollin, Lower bounds for class numbers of
real quadratic and biquadratic fields,
 Proc.  Amer.  Math. Soc. {\bf 101} (1987) 439-444.
\item[[2]] FENG Keqin, HU  Weiqun, On real quadratic function fields of Chowla
type with ideal class number one, Proc. Amer. Math. Soc.,  {\bf
127} (1999), 1301-1307.
\item[[3]] JI Guangheng, LU Hongwen, Proof of class number formula by machine, Sci. in
China,  (A){\bf 28} (1998),  193-200.
\item[[4]] S. Louboutin,  Continued fraction and real quadratic fields, J.
Number Theory, {\bf 30} (1998),  167-176.
\item[[5]] LU Hongwen, Gauss' conjectures on the quadratic number fields,
Shanghai Sci.  Tech. Pub.,  1991.
\item[[6]] ZHANG Xianke, L. C. Washington, Ideal class-groups and there
subgroups of real quadratic fields, Sci. in China,  (A){\bf 27} (1997), 522-528.
\item[[7]] E. Artin, Quadratische K\"orper im Gebiete
der h\"oheren Kongruenzen I, II, Math. Z., {\bf 19} (1924),
153-206, 207-246.
\item[[8]] WANG Kunpeng, ZHANG Xianke, The continued fractions connected with
quadratic function fields (to appear).

\end{description}
  \par \vskip 0.3cm

\end{document}